\documentclass[preprint,12pt]{article}

\usepackage{amssymb}
\usepackage{amsmath}
\usepackage{lineno}

\usepackage[normalem]{ulem}
\newcommand{\stkout}[1]{\ifmmode\text{\sout{\ensuremath{#1}}}\else\sout{#1}\fi}

\usepackage{slashbox}

\usepackage{soul}
\usepackage{lipsum}
\usepackage{amsfonts}
\usepackage{graphicx}
\usepackage{epstopdf}
\usepackage{mathtools}
\usepackage{csvsimple}
\usepackage{comment}
\usepackage{subcaption}
\usepackage{rotating}
\usepackage{arydshln}
\usepackage{url}
\usepackage{color}
\usepackage{rotating}

\usepackage{amssymb,amsmath,amsthm}
\newtheorem{theorem}{Theorem}
\newtheorem{lemma}{Lemma}
\newtheorem{corollary}{Corollary}

\usepackage{mathtools} 

\usepackage{algorithm}
\usepackage{algpseudocode}
\algrenewcommand\textproc{}

\newcommand {\rank}     {\mathop{\rm rank}\nolimits}

\newcommand {\eproof}
      {\space
        {\ \vbox{\hrule\hbox{\vrule height1.3ex\hskip0.8ex\vrule}\hrule}}
        \par}
\author{Murat Manguo\u{g}lu\footnotemark[2]
        \and Volker Mehrmann\footnotemark[3]}

\title{Multilayer approximate nullspace methods for saddle point systems}

\begin{document}

\maketitle

\renewcommand{\thefootnote}{\fnsymbol{footnote}}

\footnotetext[2]{Institut f\"{u}r Mathematik,  Technische Universit\"{a}t Berlin,  10623 Berlin,  Germany.
Present address: Department  of  Computer  Engineering,  Middle  East  Technical  University,  06800  Ankara, Turkey ({\tt manguoglu@ceng.metu.edu.tr}).}
\footnotetext[3]{Institut f\"{u}r Mathematik,  Technische Universit\"{a}t  Berlin,  10623 Berlin,  Germany ({\tt mehrmann@math.tu-berlin.de}).}

\begin{abstract}
We propose a new class of multi-layer iterative schemes for solving sparse linear systems in saddle point structure. The new scheme consist of  an iterative preconditioner that is based on the (approximate) nullspace method, combined with an iterative least squares approach and an iterative projection method. We present a theoretical analysis and demonstrate the effectiveness and robustness of the new scheme on sparse matrices from various applications.
\end{abstract}


\vskip .1truecm

\begin{abstract}
We propose a new class of multi-layer iterative schemes for solving sparse linear systems in saddle point structure. The new scheme consists of  an iterative preconditioner that is based on the (approximate) nullspace method, combined with an iterative least squares approach and an iterative projection method. We present a theoretical analysis and demonstrate the effectiveness and robustness of the new scheme on sparse matrices from various applications.
\end{abstract}

{\bf Dedicated to Daniel Szyld on the occasion of his 70th birthday.}
\section{Introduction}
We study the numerical solution of large sparse linear systems with saddle point block structure
\begin{equation}
\label{eq:general:saddle}
\mathcal W \begin{bmatrix}
x \\
y
\end{bmatrix}=\begin{bmatrix}
A & B \\
-C^T & 0
\end{bmatrix}
\begin{bmatrix}
x \\
y
\end{bmatrix}=
\begin{bmatrix}
f \\
g
\end{bmatrix},
\end{equation}
where
$A\in\mathbb{R}^{n\times n}$, 
$B,C\in\mathbb{R}^{n\times m}$ with $m\leq n$ and $\rank B=\rank C=k\leq m$. Furthermore, we assume that $A+A^T$ is positive semidefinite (denoted as $A+A^T\geq 0$).
In our theoretical analysis, we allow that $k<m$, which means that the system is singular.
In the presented numerical methods, however,  we assume that the linear system has a unique solution, which implies that  $B$ and $C$ are of rank $k = m$, since this is the property that holds in almost all of the applications.  { We refer  to \cite{BenGL05} for the necessary and sufficient conditions for the nonsingularity of $\mathcal{W}$ under various assumptions.} However, we will also discuss the treatment of the case $k<m$.

Typical applications where linear systems of the form \eqref{eq:general:saddle} arise are in the time integration of dissipative Hamiltonian systems, see \cite{GudLMS22,MehM19}, in computational fluid dynamics \cite{Bre74,elman2014ifiss}, and in constrained optimization problems \cite{GilMW21}.
In most applications, one has $C=B$
but the case $B\neq C$ also arises in practice, see the examples presented in Section~\ref{sec:testproblems}.
Very often the two components of the solution vector  represent different  quantities, such as e.g. a state vector $x$  and a Lagrange multiplier $y$, that penalizes the constraints.

We will discuss preconditioned iterative methods for the solution of \eqref{eq:general:saddle}, where typically for the preconditioner there are  two major approaches: those that use information from the underlying physical application, if they are available, i.e., one uses problem specific techniques, or those that are based on purely algebraic information, which is the approach that we study in this paper.
For detailed surveys on the iterative and direct numerical solution of problems with saddle point structure, see \cite{BenGL05,Roz18}.

In order to address many different applications, we will develop a new, purely algebraic, black-box iterative solver class for the numerical solution of large sparse  problems of the form \eqref{eq:general:saddle},
and we study  { three different subclasses:
\begin{itemize}
    \item [(i)] the  'symmetric case', where $B = C$ and $A = A^T\geq 0$, \item [(ii)] the 'structurally symmetric case', where $B=C$ but $A \not = A^T$ and $A+A^T\geq 0$, and
    \item [(iii)] the 'general case', where $B\not =C$,
    $A \not = A^T$ and $A+A^T\geq 0$.
\end{itemize}
}

A classical purely algebraic approach to the solution or the  construction of preconditioners for \eqref{eq:general:saddle} is based on the (approximate) block $LU$-factorization of the coefficient matrix, see \cite{GolW98}, that leads to a well-known scheme which then requires the approximate iterative solution of a smaller linear system where the coefficient matrix is an approximation of the Schur complement $S=C^TA^{-1}B$ which is, however, typically a dense matrix. Another common technique
is the nullspace approach, see \cite{ScoT22}, where one projects the problem onto a subproblem that is associated to the nullspace of $C^T$ which is often easily identified in practice.  While, in general for large-scale problems,  computing the exact nullspace basis is very difficult \cite{coleman1986null}, for some applications, it can be easily identified (see for example the list of applications and references given in  \cite{BenGL05}). In the complement of this nullspace, the solution is fixed by the second equation in \eqref{eq:general:saddle}, so that the problem can be reduced to the nullspace. Note that the nullspace method has recently been extended to the case of a nonzero $(2,2)$ block, see~\cite{ScoT22}.

As a modification to these well-known approaches of the nullspace method, our new method is an inner-outer iteration scheme that is based on an approximation of the classical \emph{nullspace method}, see \cite{ScoT22}, and the recently introduced \emph{shifted skew-symmetric preconditioning} \cite{ManM21}.
It employs iterative schemes for the occurring least squares problems, combined with classical projection type iterative methods, see e.g. \cite{BenCT00,SchZ07,Wat15}, as an outer iteration.
The method will be particularly interesting in the context of implicit discretization methods for linear time-varying or nonlinear  dissipative Hamiltonian differential equations, where a sequence of slightly varying systems has to be solved; see \cite{GudLMS22,MehM19}.

{Given a linear system  of the form \eqref{eq:general:saddle},
suppose that we can easily compute matrices $Z$ and $U\in\mathbb{R}^{n\times (n-m)}$ whose columns form (approximate) bases for the nullspaces of $B^T$ and $C^T$, respectively (i.e.,  $B^TZ = 0$ and $C^TU = 0$).

The second block equation  in \eqref{eq:general:saddle} is an underdetermined linear system $-C^Tx=g$ which may have infinitely many solutions.  Suppose that a particular solution $\hat{x}$ is given such that $-C^T\hat{x}=g$.
Then, with $x = \hat{x} + \bar{x}$,
we can solve the equivalent system
\begin{equation}
\begin{bmatrix}
A & B \\
-C^T & 0
\end{bmatrix}
\begin{bmatrix}
\bar{x} \\
y
\end{bmatrix}=
\begin{bmatrix}
f-A\hat{x} \\
0
\end{bmatrix}
\end{equation}
If we solve the second block row by choosing a  solution $z$ from the nullspace of $C^T $, i.e., by taking  $\bar{x}=Uz$ and  then substituting this into the  first block row, we obtain
\begin{equation}\label{eq:firstblock}
    AUz + By = f -A\hat{x}.
\end{equation}
Then, multiplying both sides of \eqref{eq:firstblock} by $Z^T$ from the left, we obtain a reduced system of the form
\begin{equation}\label{reducedsystem}
    Z^TAUz  = Z^T(f -A\hat{x})
\end{equation}
which is reflecting the two nullspaces.
Solving this reduced system for $z$, we obtain $\bar{x} = Uz$, then $x =\bar{x} +\hat{x}$, and finally  $y$ by solving the overdetermined system
\begin{equation}
    By = f -Ax.
\end{equation}

Note that to solve the over- and under-determined systems , in general, we may use a linear least squares approach via an approximate or exact sparse $QR$ decomposition , or iteratively via a Krylov subspace method.
The basic nullspace method is summarized in Algorithm~\ref{alg:general_nullspace}. %
\begin{algorithm}[H]
\begin{algorithmic}[1]
\caption{Basic nullspace method.}
\label{alg:general_nullspace}
\Require $A\in\mathbb{R}^{n\times n}$, $B,C\in\mathbb{R}^{n\times m}$, $f\in \mathbb{R}^n$ and $g\in \mathbb{R}^m$ with $m<n$.
\Ensure $x\in \mathbb{R}^n$ and $y\in \mathbb{R}^m$
\State{Construct $Z, U\in\mathbb{R}^{n\times (n-m)}$ whose columns form  bases for the nullspaces of $B^T$ and $C^T$, respectively.}
\State{Find a particular solution, $\hat{x}\in\mathbb{R}^n$, of  $-C^Tx=g$.
}
\State{Solve   $Z^TAUz  = Z^T(f -A\hat{x})$.}
\State{Set $x =\hat{x}+ Uz$.}
\State{Solve $By = f - Ax$.}
\end{algorithmic}
\end{algorithm}
Note that in the symmetric and structurally symmetric case, where $C=B$, the method simplifies, since one needs to compute only one nullspace basis $Z$, and lines 3 and 4 change to $Z^TAZz=Z^T(f-A\hat{x})$ and $x = \hat{x}+Zz$, respectively. We call $\mathcal{N}=Z^TAZ$ and $\mathcal{N}=Z^TAU$ the projected nullspace matrices, for the symmetric/structurally symmetric and general cases, respectively.  { We note that $\mathcal{N}$ is also known as the {\em reduced Hessian} matrix in optimization (see for example \cite{nocedal2006numerical})}.
In Section~\ref{sec:theory} we analyze the properties of the symmetric part of the projected (approximate) nullspace matrix, i.e.,  $\mathcal{N}^s=\frac{Z^TAU+(Z^TAU)^T}{2}$ and show under which conditions it is positive definite.

The nullspace method is exact only if the nullspace bases are determined exactly, and the reduced under- as well as over-determined systems are solved  exactly.
Since we cannot expect this to happen, we propose to compute all of these solutions approximately and use an outer iterative projection based  Krylov subspace  method to correct the approximations. In this way, the approximate general nullspace method becomes a method to apply a preconditioner. This approximate approach also  addresses a  major disadvantage of the exact nullspace method, that the matrices $Z$ and $U$ are usually not sparse, even if $B$ and $C$ are sparse. In the approximate nullspace approach, we can employ  sparse approximate nullspace bases,
see~\cite{ScoT22} for a detailed survey of the available methods.

The remainder of the paper is organized as follows. In Section \ref{sec:theory} we present the underlying theoretical foundations of the new schemes. Then, in Section \ref{sec:ourmethod} we describe the new schemes. The description of the test problems, the experimental setup and the results are given in Sections \ref{sec:testproblems}, \ref{sec:framework}, and \ref{sec:results}, respectively. We conclude and discuss future work in Section \ref{sec:conclusions}.

\section{Theoretical Analysis}\label{sec:theory}
We will make frequent use of the following important results, see \cite{AchAM21,GudLMS22}, on matrices with positive semidefinite symmetric part, presented in a form needed in our analysis.
\begin{lemma}\label{lem:sts}
Consider $A=J+R\in\mathbb R^{n\times n}$, where $0\leq R=R^T $ and $0\neq J=-J^T$. Then there exist a real orthogonal matrix $Q\in\mathbb R^{n
\times n}$, and integers $n_1 \geq n_2 \geq \cdots \geq n_{r-1} >0$ and $n_r,n_{r+1} \geq 0$, such that
\begin{equation}\label{eq:stcform}
Q^TRQ= \begin{bmatrix} R_{11} & 0 \\ 0 & 0\end{bmatrix},\
Q^TJ Q=
\begin{bmatrix}
J_{11} & J_{12} &  & & 0 &0\\
J_{21} & J_{22} &  \ddots &  & 0 &0\\
&  \ddots & \ddots  & J_{r-2,r-1} &\vdots &\vdots \\
& &  J_{r-1,r-2}& J_{r-1,r-1} & 0 & 0\\
0 & \cdots & \cdots  & 0 & J_{r,r} &0\\
0& \cdots & \cdots &0 & 0 & 0
\end{bmatrix} ,
\end{equation}
where $0<R_{11}=R_{11}^T \in\mathbb R^{n_1\times n_1}$ , $J_{ii}=-J_{ii}^T\in\mathbb R^{n_i\times n_i}$ for $i=1,\dots,r$, $J_{rr}$ invertible,
and $J_{i,i-1}=-J_{i-1,i}^T=[\Sigma_{i,i-1}\;0]\in\mathbb R^{n_{i}\times n_{i-1}}$ with $ \Sigma_{i,i-1}$
being nonsingular for $i=2,\dots,r-1$.
\end{lemma}
\proof
The proof follows directly from the staircase form of \cite{AchAM21} by a further full rank decomposition of the bottom right block.
\eproof

For a system transformed to the form \eqref{eq:stcform}, one can perform a block $LDU$ factorization, see \cite{GudLMS22}.
\begin{lemma}\label{lem:schur}
In the notation of Lemma~\ref{lem:sts}, the matrix $Q^TAQ$ can be transformed via Schur complement reduction into the block diagonal form
\begin{equation}\begin{bmatrix}
\widehat{A}_{11} & &  & & & \\
 & \widehat{A}_{22} &   &  & & \\
&  & \ddots  &  &  &\\
& &  & \widehat{A}_{r-1} & &\\
 & &  & & \widehat {A}_{r,r}&\\
 &&&&&0
\end{bmatrix},
\label{eq:schurform}
\end{equation}
where for $i=1,\ldots,r-1$, $\widehat{A}_{ii}=\widehat R_{ii}+\widehat J_{ii}$ with $\widehat R_{ii}>0$ and $\widehat A_{r,r}=-\widehat A_{r,r}^T$ is invertible. Note that the last two diagonal blocks may not be present, and in this case the symmetric part of \eqref{eq:schurform} is positive definite.
\end{lemma}
\proof See \cite{GudLMS22}.
\eproof
Based on these two lemmas we then have the following  result.
\begin{theorem}\label{thm:nspace}
Consider a symmetric, structurally symmetric or general linear system of the form~\eqref{eq:general:saddle}, with $A+A^T\geq 0$, 
and $\rank B=\rank C^T= k\leq m$. Then there exist real orthogonal matrices $Q, V$ such that
\begin{equation}\label{eq:blockform}
    Q^T\mathcal{W}V=
    \begin{bmatrix}
    A_{11} & A_{12} & A_{13} &A_{14}&  B_{11} &0 \\
    A_{21} & A_{22} &0 &0&B_{21} &0\\
    -A_{13}^T & 0&A_{33} & 0 & B_{31} &0\\
    -A_{14}^T & 0 & 0 & 0 & B_{41} &0 \\
    -C_{11}^T & 0 & 0 & 0 & 0 & 0\\
    0 & 0 & 0 & 0 & 0 & 0
    \end{bmatrix},\
V^T   \begin{bmatrix}x\\ y\end{bmatrix}=
\begin{bmatrix} x_1\\ x_2\\ x_3\\ x_4\\ y_1 \\ y_2
\end{bmatrix},\
Q^T   \begin{bmatrix}f\\ g \end{bmatrix}=
\begin{bmatrix} f_1\\ f_2\\
f_3\\ f_4\\ g_1 \\ g_2
\end{bmatrix},
\end{equation}
with $A_{13}=-A_{31}^T$, the matrices $A_{33}=-A_{33}^T\in \mathbb R^{n_3 \times n_3}$,
$B_{11}, C_{11}\in \mathbb R^{n_1 \times n_1}$, $A_{22}\in \mathbb R^{n_2\times n_2}$  are invertible, and $A_{22}$ is
of the form~\eqref{eq:stcform}.

Moreover, in the symmetric and structurally symmetric case $B_{21}, B_{31}, B_{41}=0$.
\end{theorem}
\proof
    The proof is constructive by a sequence of equivalence and congruence transformations with real orthogonal transformations. It can be directly implemented in a numerically stable method.

   Step 1. Perform a singular value decomposition
   \[
   Q_2^T C^T B V_2= \begin{bmatrix} K_{1} & 0 \\ 0 & 0 \end{bmatrix},
   \]
   where $K_1\in \mathbb R^{n_1 \times n_1}$ is square and nonsingular and $Q_2,V_2$ are orthogonal.

Step 2. Perform a $QR$ decomposition with $Q_1$ orthogonal such  that
   \[
   (Q_2^T C^T) Q_1 = \begin{bmatrix} \hat C_{11}^T & 0\\ \hat C_{12}^T  & 0
   \end{bmatrix},\ Q_1^T B V_2 = \begin{bmatrix} \hat B_{11} & \hat B_{12} \\ \hat B_{21} & \hat B_{22}
   \end{bmatrix}.
   \]
Then it follows that  $\hat C_{11}^T\hat B_{11}=K_1$ is nonsingular and hence $\hat C_{11},\hat B_{11}$ are nonsingular, as well as $\hat B_{12}=0$, $\hat C_{12}=0$. Then the assumption that
$\rank C=\rank B=k$ implies that $\hat B_{22}=0$.
Let
\[Q_1^T A Q_1= \begin{bmatrix} \hat A_{11} & \hat A_{12} \\ \hat A_{21} & \hat A_{22}
   \end{bmatrix},
\]
and apply Lemma~\ref{lem:sts} to determine an orthogonal matrix $Q_{22}$ that generates a decomposition of $\hat A_{22}$ in the form \eqref{eq:stcform}. Then \eqref{eq:blockform} is obtained with
\[ Q= \begin{bmatrix} I_{n} & 0 \\ 0& Q_{2}
\end{bmatrix} \begin{bmatrix} I_{k} & 0 &0\\ 0& Q_{22} & 0 \\ 0 & 0 & I_{m-k}
\end{bmatrix},\  V= \begin{bmatrix} I_{n} & 0 \\ 0& V_2
\end{bmatrix} \begin{bmatrix} I_{k} & 0 &0\\ 0& Q_{22} &0 \\ 0 & 0 & I_{m-k}
\end{bmatrix}.
\]
\eproof
 {Note that the transformation in Theorem~\ref{thm:nspace} is an orthogonal congruence on the $A$ block and hence the property that $A+A^T\geq 0$ is preserved under the transformation. Note further that if $B=C$ then $Q_2=V_2$ and thus $Q=V$ and  $B_{21}, B_{31}, B_{41}=0$.}

Theorem~\ref{thm:nspace} directly implies the following properties of the solution for a system of the form \eqref{eq:blockform}.
The vector $y_2$ can be chosen arbitrarily, so for the solution to be unique, we need that the last block row and column is void.
This is the case if and only if $k=m$.

In the symmetric and structurally symmetric case
we also have that $x_1= -C_{11}^{-T}g_1$,
$x_3=A_{33}^{-1} (A_{31} C_{11}^{-T} g_1+f_3)$,
$x_2= A_{22}^{-1} (-A_{21} x_1 +f_2)$, and $x_4$ is arbitrary.
So if $\mathcal{W}$ is invertible, then the forth row and column must be void as well.

We have the following consequences.
\begin{corollary}\label{cor:hypo}
~
Consider a linear system of the form ~\eqref{eq:general:saddle} with $\mathcal W$ invertible.
i) For the symmetric case, i.e.,  if $A=A^T\geq 0$,  $B=C$,  then the nullspace matrix satisfies $Z^TAZ>0$. \\
ii)  For the structurally symmetric  case, i.e., if  $A\neq A^T$,  but $A+A^T\geq 0$, $B=C$, and   $\mathcal W$  is invertible, then there exists a real orthogonal matrix $Q$ such that  $Q^TZ^TAZQ=\mbox{\rm diag}(A_{22} ,A_{33})$, with $A_{33}$ skew-symmetric and invertible. If the matrix $A$
 { can be transformed to the form \eqref{eq:stcform} with the last two rows and columns void},
then the block $A_{33}$ is void and for $A_{22}$ there exists a block $LDU$ factorization as in \eqref{eq:schurform}, where the last block is void and the symmetric part is positive definite. \\
iii)  For the general  case, i.e., if  $A\neq A^T$,  but $A+A^T\geq 0$, $B\neq C$, if the matrix $A$  { can be transformed to the form \eqref{eq:stcform} with the last two rows and columns void}, then the block $A_{33}$ is void and for $A_{22}$ there exists a block $LDU$ factorization as in \eqref{eq:schurform}, where the last block is void and the symmetric part is positive definite.
\end{corollary}

Note that in our algorithmic approaches, the condition in Corollary~\ref{cor:hypo}
that $A+A^T\geq 0$ can be relaxed that it holds only in the nullspace $Z$.

The construction of the projection matrices $Z,U$ with orthonormal columns will typically lead to full matrices. So for the practical implementation we will use sparse approximate versions of the associated decompositions. In this way we will not obtain exact bases for the nullspaces, but sparse approximations $\tilde Z,\tilde U$. For the symmetric and structurally symmetric case, the following result shows
that if the error $\|Z-\tilde Z\|$ can be bounded by  a small quantity, then the approximate nullspace matrix $\tilde Z^T A\tilde Z$ will still have a positive definite symmetric part.

\begin{lemma}\label{lem:perturbed}
Under the assumptions of Corollary~\ref{cor:hypo}, let in the symmetric and structurally symmetric case $\tilde Z$ be a (sparse) approximation of $Z$ (with columns that are not exactly orthonormal) such that $\|\tilde Z-Z\|_2\leq \epsilon \|Z\|_2= \epsilon$, with a sufficiently small $\epsilon $. Then $\tilde Z^TA \tilde Z$ still has a positive definite symmetric part.
    \end{lemma}
\proof
Let $\Delta Z=\tilde Z-Z$ and use the fact that the eigenvalues $\tilde Z^T A \tilde Z$ are continuous functions in the entries of $\Delta Z$. Since the eigenvalues of $Z^TAZ$ have  positive real parts, it then follows that for sufficiently small $\epsilon$ the eigenvalues of $\tilde Z^T A \tilde Z$ have positive real parts as well.
\eproof
Note that in the general case it is difficult to make a statement about the solvability and the symmetric part of the nullspace matrix, although, as we demonstrate in Section~\ref{sec:results}, in the discussed applications it is typically positive definite.

\section{A new multi-layer iterative scheme}\label{sec:ourmethod}

The main motivation and core idea of our new iterative multilayer scheme  is the availability of efficient and robust Krylov subspace methods for shifted skew-symmetric~\cite{GudLMS22,IdeV07,IdeV23,Jia07,Rap78,Wid78}  as well as symmetric and positive definite~\cite{GolV96,hestenes1952methods,Saa93}} systems. These methods have short recurrences and satisfy optimality properties. We make use of these properties and propose a robust scheme based on preconditioning via the approximate nullspace method.

While for the solution of saddle point problems in the symmetric and structurally symmetric case, nullspace methods
have been studied extensively \cite{Roz18,ScoT22}, the general case has not received much attention, due to the associated computational and memory cost, which involves  two different nullspace bases. Typically, furthermore, for the usual nullspace approaches it is required that the nullspace dimension is small.
In contrast to this, by allowing sparse approximations, we show how these issues can be alleviated in our new technique.

We determine sparse approximate bases for the nullspaces of $B^T$ (or $C^T$) by the right oblique conjugation scheme~\cite{ScoT22}. This approach
can be viewed as a 'partial' biconjugation~\cite{ChuFG95} and direct projection method~\cite{BenM95,ScoT22}, such that given $B\in \mathbb{R}^{n\times m}$
(with $\rank(B^T)=k\leq m$), assuming that no pivoting is required,
\begin{equation}
    B^TV = L,
\end{equation}
where $V\in \mathbb{R}^{n\times n}$ is invertible and
$L\in \mathbb{R}^{k\times n}$ is a lower trapezoidal matrix. Then, the rightmost $n-k$ columns of $V$ (i.e., $Z=V(:,k+1:n)$  in MATLAB notation) form a basis for the nullspace. Since in our case, we are not interested in $L$, it is not stored. However, since $B$ is sparse, pivoting is usually required to avoid division by zero or small values in modulus. The construction for the nullspace of $C$ is analogous.

A pseudocode for the sparse approximate right oblique conjugation with pivoting to obtain a sparse approximate nullspace basis is presented in Algorithm \ref{alg:oblique}. Hereafter, we refer to this procedure as \emph{SAROC}.  It produces a sparse basis, since we assume that computed coefficients that are smaller than a given threshold and entries of the matrix for which we compute the nullspace that are lower than a given tolerance in absolute value are set to zero, as shown in lines 8 and 9, respectively.
As a consequence, we  obtain \emph{approximate} projected nullspace matrices, $\widetilde{\mathcal{N}}=\widetilde{Z}^TA\widetilde{Z}$ and $\widetilde{\mathcal{N}}=\widetilde{Z}^T A\widetilde{U}$, respectively, for the  { symmetric/}structurally symmetric and general case. Assuming dense operations, the worst-case cost of SAROC is $\mathcal{O}(n^2)$  arithmetic operations per iteration, giving a total cost of $\mathcal{O}(kn^2)$.  However, since the  matrix $V$ is initialized as the identity matrix and the algorithm is designed to maintain sparsity via the drop tolerance ($\tau$) and the threshold  ($\rho$)-controlled vector updates, the cost is reduced significantly. Let $s_v$ and $s_b$ denote the average number of nonzeros per column of $V$ and $B$, respectively. Inner products (line 3) cost  $\mathcal{O}(sn)$ per iteration and vector updates (line 9) cost $\mathcal{O}(s_v)$ per update, where $s=\min(s_v,s_b)$. Since the threshold $\rho$ limits the updates  (worst-case $n-i$ per iteration), the total cost would be approximately $\mathcal{O}(ksn)$ where $s \ll n$. Furthermore, since the algorithm performs column pivoting (lines 4-6), the matrix $V$ must be stored in column-major order to avoid  cache misses.

\begin{algorithm}[H]
\begin{algorithmic}[1]
\caption{Sparse Approximate Right Oblique Conjugation (SAROC) for computing an approximate nullspace basis}
\label{alg:oblique}
\Require $B\in \mathbb{R}^{n\times m}$ with  $\rank(B^T)=k\leq m$, $\rho$ (threshold) and $\tau$ (drop tolerance)   \Comment{$b_{i}$: is the $i^{th}$ column of $B$}
\Ensure $\widetilde{Z}\in\mathbb{R}^{n\times (n-k)}$  \Comment{$z_{i}$: is the $i^{th}$ column of $\widetilde{Z}$}
\State{$v_i =e_i$ for $(1\leq i\leq n)$}
\For{$i = 1,\ldots,k$}
\State{$\sigma  =  b_i^T[v_i, v_{i+1},\ldots, v_n]$}  \Comment{compute coefficient vector $\sigma \in \mathbb{R}^{n-i+1}$}
\State{$l^* = \arg\max\limits_{l}{(|\sigma_l|)}$}\Comment {find  a good pivot}
\State{$v_{i} \leftrightarrow v_{i+l^*-1}$} \Comment{Swap the corresponding columns of $V$}
\State{$\sigma_{1} \leftrightarrow \sigma_{l^*}$} \Comment{Swap the corresponding coefficients}
\For{$j = i+1,\ldots,n$}
\If{$|\sigma_{j-i+1}/\sigma_{1}|>\rho$}
\State{$v_j \xleftarrow{\text{$|v_j(\cdot)|\ge\tau||v_j||_2$}} v_ j - (\sigma_{j-i+1}/\sigma_{1})v_i$}
\EndIf
\EndFor
\EndFor
\State{$\widetilde{Z} = [v_{k+1},v_{k+2}, \ldots, v_{n}]$}
\end{algorithmic}
\end{algorithm}
The next procedure common to all cases is to compute the Factorized (Cholesky) Sparse Approximate Inverse (FSAI) of the symmetric part $\widetilde{\mathcal N}^s$ of the approximate projected nullspace matrix. This step is required to reveal the shifted skew-symmetric system in Step 3 of Algorithm~\ref{alg:general_nullspace}.

Rather than computing the classical Cholesky factorization of the symmetric part $\widetilde{\mathcal N}^s$ of the approximate projected nullspace matrix,  to avoid explicitly forming this matrix,  we opt for a factorized approximate inverse based Cholesky factorization. If $\widetilde{\mathcal N}^s$ is not positive definite, then the Cholesky factorization will detect this. If $\widetilde{\mathcal N}^s$ is indefinite (with a few negative eigenvalues), it would still be possible to obtain a shifted-skew symmetric system, as shown in~\cite{ManM21}, but with an additional cost.  { Alternatively, we could just apply a diagonal shift  $\widetilde{\mathcal N}^s + \alpha I$ but  with an additional cost of finding a 'good' $\alpha$, which can be obtained, for example, by estimating the smallest eigenvalue.}

We obtain the FSAI of $\widetilde{\mathcal{N}}^s$ using the method in~\cite{BenCT00}. Our implementation, presented in Algorithm \ref{alg:fsai},  is implicit in the sense that $\widetilde{\mathcal{N}}^s$ is not formed explicitly, since it can be relatively dense even if $\widetilde{Z}$, $\widetilde{U}$, and $A$ are sparse. Instead,  its columns are extracted, as needed in line 3. This can be easily achieved by implicit matrix-vector multiplications, $n^s_i= \widetilde{\mathcal{N}}^s e_i $. Similar to the sparse approximate right oblique conjugation, we replace computed  coefficients that are smaller than a given threshold and entries of the matrix that is factorized that are smaller than a given tolerance in absolute value  by zero, see lines 5 and 6, respectively.   FSAI  produces $\widetilde{W}\in\mathbb{R}^{(n-k)\times (n-k)}$ such that $\widetilde{W}^T\widetilde{\mathcal{N}}^s\widetilde{W} \approx I$. Furthermore, $\widetilde{\mathcal{N}}^s$ and also the skew-symmetric part  $\widetilde{\mathcal{N}}^j$, or any other matrices resulting from further projections on them, do not need to be formed explicitly,  since matrix-vector multiplications with them can be computed implicitly.

Due to approximations and thresholding, Algorithm~\ref{alg:general_nullspace} is no longer a direct factorization,  but rather a preconditioner, and hence we need  multiple layers of iterative methods. We use Krylov subspace methods throughout our new schemes, namely: the restarted \emph{flexible Generalized Minimal Residual method}  (fGMRES)~\cite{saad1993flexible} for general systems, where the preconditioner is  another iterative solver,  the \emph{Minimal Residual method for Shifted Skew-Symmetric Systems} (MRS) for shifted skew-symmetric systems~\cite{ManM21},  the \emph{Conjugate Gradient method} (CG) for symmetric and positive definite systems~\cite{GolV96}, and the \emph{Least squares $QR$ method} (LSQR)~\cite{paige1982lsqr} for over-determined and under-determined systems.

In the following subsections,  we describe each of the resulting schemes in detail.
\begin{algorithm}[H]
\begin{algorithmic}[1]
\caption{FSAI (A-Conjugation)}
\label{alg:fsai}
\Require $\widetilde{\mathcal{N}}^s\in\mathbb{R}^{(n-k)\times(n-k)}$, $\rho$ (threshold) and $\tau$ (drop tolerance)
\Ensure $\widetilde{W}\in\mathbb{R}^{n\times (n-k)}$
\State{$w_i =e_i$ for $(1\leq i\leq n-k)$}
\For{$j = 1,\ldots,n-k$}
\State{$[\sigma_j, \sigma_{j+1},\ldots, \sigma_{n-k}]^T = [w_j,w_{j+1},\ldots, w_{n-k}]^T n^s_{j}$} \newline\Comment{$n^s_{j}$ is the $j^{th}$ column of $\widetilde{\mathcal{N}}^s$  and $w_{i}$ is the $i^{th}$ column of $\widetilde{W}$}
\For{$i = j+1,\ldots,n-k$}
\If{$|\sigma_{i}/\sigma_{j}|>\rho$}
\State{$w_i \xleftarrow{\text{$|w_i(\cdot)|\ge\tau||w_i||_2$}} w_i- (\sigma_i/\sigma_j)w_j$}
\EndIf
\EndFor
\EndFor
\State{$\sigma_i = 1/\sqrt{\sigma_i}$ for $(1\leq i\leq n-k)$}
\State{$\widetilde{W} = \widetilde{W}P$ where $P$ is a diagonal matrix with entries $P_{(i,i)} = \sigma_i$}
\end{algorithmic}
\end{algorithm}

\subsection{Symmetric case}\label{sec:i)}
In the symmetric case, we have $B=C$ and $A=A^T$, and we compute a sparse approximate nullspace basis, $\widetilde{Z}$,  of $B^T$.
We then assume, see Section~\ref{sec:theory} for conditions, that ${\mathcal{N}}^s=\widetilde{Z}^TA\widetilde{Z}$ is positive definite. Since the new scheme involves sparse approximations and is just used as a preconditioner, in the following numerical experiments, we assume  that $\widetilde{\mathcal{N}}^s$ is positive definite, even if $A$ is indefinite or semidefinite. This allows us to cover a larger number applications than in our initial assumption. The pseudocode of the new scheme for the symmetric case is given in Algorithm~\ref{alg:symmetric_nullspace}.

Optionally,
we can further $M$-orthogonalize the columns of $\widetilde{Z}$ with respect to $M=A$ using the modified Gram-Schmidt (MGS) procedure given in Algorithm~\ref{alg:Morth-MGS}. Since it may be too expensive to orthogonalize against all vectors, we also adopt a 'partial' orthogonalization against a certain number of previous vectors, which we refer to as the window size $w$ (used in line 4) and we employ numerical dropping with a threshold $\tau$ (used in line 7) to ensure sparsity. Although $M$-orthogonalization can be achieved during the oblique conjugation process in Algorithm~\ref{alg:oblique} for the last $(n-k)$ columns of $V$, we prefer to keep it as a separate stage in our new scheme, since it requires the matrix $M=A$ which may change when solving multiple linear systems with different  right-hand sides even though $B$ remains the same.

Note that FSAI in step 2 is already computing an $M$-orthgonolization of $\widetilde{Z}$, i.e., the columns of $\widetilde{Z}\widetilde{W}$ are expected to be orthogonal with respect to $M=A$,  since $\widetilde{W}^T\widetilde{Z}^TA\widetilde{Z}\widetilde{W}\approx I$ and $B^T\widetilde{Z}\widetilde{W}\approx 0$. Nevertheless, since both are approximations, we keep both versions and provide MGS based $M$-orthogonlization as an option.

\begin{algorithm}[ht]
\begin{algorithmic}[1]
\caption{MGS $M$-orthogonalization.}
\label{alg:Morth-MGS}
\Require $Z\in\mathbb{R}^{n\times (n-k)}$, $M\in\mathbb{R}^{n\times n}$  \Comment{$z_{i}$: is the $i^{th}$ column of $Z$}
\Ensure $\bar{Z}\in\mathbb{R}^{n\times (n-k)}$  \Comment{$\bar{z}_{i}$: is the $i^{th}$ column of $\bar{Z}$}
\State{$\bar{z}_1 =z_1/\sqrt{z_1^TMz_1}$}
\For{$i = 2,\ldots,n-k$}
\State{$\bar{z}_i = z_i$}
\For{$j = \max(i-w,1),\ldots,i-1$}
\State{$\bar{z}_i = \bar{z}_i - (\bar{z}_j^TM\bar{z}_i)\bar{z}_j$}
\EndFor
\State{$\bar{z}_i(\cdot) \xleftarrow{\text{$|\bar{z}_i(\cdot)|<\tau||\bar{z}_i||_2$}} 0$}
\State{$\bar{z}_i =\bar{z}_i/\sqrt{\bar{z}_i^TM\bar{z}_i}$}
\EndFor
\end{algorithmic}
\end{algorithm}
Note that our assumption that $Z^TAZ>0$ implies that in line 8 of Algorithm~\ref{alg:Morth-MGS} we have ${\bar z}_i^TM {\bar z}_i>0$.

We then have the approximate nullspace method for the symmetric case given in Algorithm~\ref{alg:symmetric_nullspace}.

\begin{algorithm}[H]
\begin{algorithmic}[1]
\caption{New Scheme (symmetric case)}
\label{alg:symmetric_nullspace}
\Require $A\in\mathbb{R}^{n\times n}$ with $A=A^T$, $B\in\mathbb{R}^{n\times m}$, $f\in \mathbb{R}^n$ and $g\in \mathbb{R}^m$ with $m<n$
\Ensure $x\in \mathbb{R}^n$ and $y\in \mathbb{R}^m$
\State{Construct $\widetilde{Z}\in\mathbb{R}^{n\times (n-m)}$ whose columns forms a sparse approximate basis for the nullspace of $B^T$ via SAROC}
\State{Determine the factorized sparse approximate inverse (FSAI)  of $\widetilde{\mathcal{N}}^s=\widetilde{Z}^TA\widetilde{Z}$ such that $\widetilde{W}^T\widetilde{\mathcal{N}}^s\widetilde{W}\approx I$ }
\State{Solve {\tiny $\begin{bmatrix}
A & B \\
-B^T & 0
\end{bmatrix}
\begin{bmatrix}
x \\
y
\end{bmatrix}=
\begin{bmatrix}
f \\
g
\end{bmatrix}$} via fGMRES, where the systems that involve the preconditioner {\tiny M
$\begin{bmatrix}
z_1 \\
z_2
\end{bmatrix}=
\begin{bmatrix}
t_1 \\
t_2
\end{bmatrix}$} are solved via the nullspace method as follows:}
\begin{itemize}
    \item Find a particular solution, ,$\hat{z}_1\in\mathbb{R}^n$, of  $-B^Tz_1=t_2$ via  LSQR
    \item Solve  $\widetilde{W}^T\widetilde{Z}^TA\widetilde{Z}\widetilde{W}(\widetilde{W}^{-1}u)  = \widetilde{W}^T\widetilde{Z}^T(t_1 -A\hat{z_1})$ via CG
    \item Recover $z_1 =\hat{z_1}+ \widetilde{Z}u$
    \item Solve $Bz_2 = t_1 - Az_1$ via LSQR
\end{itemize}
\end{algorithmic}
\end{algorithm}

\subsection{Structurally symmetric case}\label{sec:generalized}
In the structurally symmetric case, we assume $B=C$  and that
$A\ne A^T$.
We compute a sparse approximate nullspace basis, $\widetilde{Z}$,  of $B^T$. We again assume that $\widetilde{\mathcal{N}}^s=\widetilde{Z}^T(\frac{A+A^T}{2})\widetilde{Z}>0$, see Section~\ref{sec:theory}. The pseudocode is given in Algorithm~\ref{alg:generalized_nullspace}. We perform matrix vector multiplications with $\widetilde{\mathcal{N}}^s$ and $\widetilde{\mathcal{N}}^j$ implicitly,  while forming the symmetric and skew-symmetric parts of $A$ explicitly.

 { Similar to the symmetric case,} optionally, we can further $M$-orthogonalize the columns of $\tilde{Z}$ with respect to the symmetric part of $A$ via the modified Gram-Schmidt (MGS) procedure given in Algorithm~\ref{alg:Morth-MGS}.

\begin{algorithm}[H]
\begin{algorithmic}[1]
\caption{New Scheme (structurally symmetric case)}
\label{alg:generalized_nullspace}
\Require $A \in\mathbb{R}^{n\times n}$ with $A\ne A^T$, $B\in\mathbb{R}^{n\times m}$, $f\in \mathbb{R}^n$ and $g\in \mathbb{R}^m$ with $m<n$
\Ensure $x\in \mathbb{R}^n$ and $y\in \mathbb{R}^m$
\State{Construct $\widetilde{Z}\in\mathbb{R}^{n\times (n-m)}$ whose columns form a sparse approximate basis for the nullspace of $B^T$ via SAROC}
\State{Determine the factorized sparse approximate inverse (FSAI)  of  $\widetilde{\mathcal{N}}^s$ = $\widetilde{Z}^T(\frac{A+A^T}{2})\widetilde{Z}$ such that $\widetilde{W}^T\widetilde{\mathcal{N}}^s\widetilde{W}\approx I$ }
\State{Solve {\tiny $\begin{bmatrix}
A & B \\
-B^T & 0
\end{bmatrix}
\begin{bmatrix}
x \\
y
\end{bmatrix}=
\begin{bmatrix}
f \\
g
\end{bmatrix}$} via fGMRES, where the systems that involve the preconditioner {\tiny M
$\begin{bmatrix}
z_1 \\
z_2
\end{bmatrix}=
\begin{bmatrix}
t_1 \\
t_2
\end{bmatrix}$} are solved via the nullspace method as follows:}
\begin{itemize}
    \item Find a particular solution, $\hat{z}_1\in\mathbb{R}^n$, of  $-B^Tz_1=t_2$ via LSQR
    \item Solve  $\widetilde{W}^T\widetilde{Z}^TA\widetilde{Z}\widetilde{W}(\widetilde{W}^{-1}u)  = \widetilde{W}^T\widetilde{Z}^T(t_1 -A\hat{z_1})$ via fGMRES  with the preconditioner $N= I + \widetilde{W}^T\mathcal{N}^J\widetilde{W}$
    (systems involving the preconditioner are solved via MRS),  where $\mathcal{N}^J=\widetilde{Z}^T(\frac{A-A^T}{2})\widetilde{Z}$
    \item Recover $z_1 =\hat{z_1}+ \widetilde{Z}u$
    \item Solve $Bz_2 = t_1 - Az_1$ via LSQR
\end{itemize}
\end{algorithmic}
\end{algorithm}

\subsection{General case} \label{sec:general}

In the general case, we allow $B\neq C$ and $A\neq A^T$ and we compute  sparse approximate nullspace bases, $\widetilde{Z}$ and $\widetilde{U}$  of $B^T$ and $C^T$, respectively. Here we again assume that the symmetric part $\widetilde{\mathcal{N}}^s=(\widetilde{Z}^TA\widetilde{U}+\widetilde{U}^TA^T\widetilde{Z})/2$, of $\widetilde{Z}^TA\widetilde{U}$, is positive definite. The pseudocode is given in Algorithm~\ref{alg:general_nullspace2}.  Similar to the structurally symmetric case, in the general case, we also perform the vector multiplications with  $\widetilde{\mathcal{N}}^s$ and $\widetilde{\mathcal{N}}^j$ implicitly, but  this requires matrix vector multiplications involving both $A$ and $A^T$.

Unlike in the symmetric and structurally symmetric case, rather than $M$-orthogonali\-zation, one needs to biorthogonalize   $\widetilde{Z}$ and $\widetilde{U}$  with respect to $\widetilde{\mathcal{N}}^s$. One way to achieve this is by computing the FSAI of $\widetilde{\mathcal{N}}^s$. In other words, $\widetilde{Z}\widetilde{W}$ and $\widetilde{U}\widetilde{W}$
are already bi-orthogonalized with respect to $\widetilde{\mathcal{N}}^s$. We do not perform any further bi-orthgonalization step in the general case.

\begin{algorithm}[H]
\begin{algorithmic}[1]
\caption{New Scheme (general case)}
\label{alg:general_nullspace2}
\Require $A\in\mathbb{R}^{n\times n}$ with $A\ne A^T$, $B, C\in\mathbb{R}^{n\times m}$ with $B\ne C$, $f\in \mathbb{R}^n$ and $g\in \mathbb{R}^m$ with $m<n$
\Ensure $x\in \mathbb{R}^n$ and $y\in \mathbb{R}^m$
\State{Construct $\widetilde{Z}$ and $\widetilde{U}\in\mathbb{R}^{n\times (n-m)}$ whose columns forms sparse approximate bases for the nullspaces of $B^T$ and $C^T$, respectively,  via SAROC}
\State{Determine the factorized sparse approximate inverse (FSAI)  of $\widetilde{\mathcal{N}}^s=(\widetilde{Z}^TA\widetilde{U}+\widetilde{U}^TA^T\widetilde{Z})/2$ such that $\widetilde{W}^T\mathcal{N}^s\widetilde{W}\approx I$
}
\State{Solve {\tiny $\begin{bmatrix}
A & B \\
-C^T & 0
\end{bmatrix}
\begin{bmatrix}
x \\
y
\end{bmatrix}=
\begin{bmatrix}
f \\
g
\end{bmatrix}$} via fGMRES, where the systems that involve the preconditioner {\tiny M
$\begin{bmatrix}
z_1 \\
z_2
\end{bmatrix}=
\begin{bmatrix}
t_1 \\
t_2
\end{bmatrix}$} are solved via the nullspace method as follows:}
\begin{itemize}
  \item Find a particular solution, ,$\hat{z}_1\in\mathbb{R}^n$, of  $-C^Tz_1=t_2$ via LSQR
    \item Solve  $\widetilde{W}^T\widetilde{Z}^TA\widetilde{U}\widetilde{W}(\widetilde{W}^{-1}u)  = \widetilde{W}^T\widetilde{Z}^T(t_1 -A\hat{z_1})$ via fGMRES  with a preconditioner $N= I + \widetilde{W}^T \widetilde{\mathcal{N}}^J\widetilde{W}$ (systems involving the preconditioner are solved via MRS) where $\widetilde{\mathcal{N}}^J=(\widetilde{Z}^TA\widetilde{U}-\widetilde{U}^TA^T\widetilde{Z})/2$
    \item Recover $z_1 =\hat{z_1}+ \widetilde{U}u$
    \item Solve $Bz_2 = t_1 - Az_1$ via LSQR
\end{itemize}
\end{algorithmic}
\end{algorithm}

\section{Test problems}\label{sec:testproblems}

To evaluate the  new solver class,  we use a set of test matrices for each problem class.

 { For the symmetric case,  these are eight test problems obtained from the SuiteSparse Matrix Collection~\cite{davis2011university} (given in Table~\ref{tab:symmetric:matrices}).  These problems are  part of an optimal control collection (Vehical Dynamics and Optimization Lab - VDOL group) and  arise in various applications.  Specifically,  reorientation, dynamic soaring, Goddart rocket, tumor angiogenesis problems are described in~\cite{betts2010practical}, \cite{zhao2004optimal}, \cite{goddard1920method} and  \cite{ledzewicz2008analysis}, respectively. In Table~\ref{tab:symmetric:matrices}, we provide the smallest and largest eigenvalues of $A$. For most of the test problems in this set, $A$ is indefinite and for two problems $A$ is  positive definite.}

For the structurally symmetric case,  eight matrices are used (given in Table~\ref{tab:generalized:matrices}). The first three, namely rajat04, rajat14 and fpga\_trans\_01, are  from the SuiteSparse Matrix Collection.   The remaining five test matrices are generated from the IFISS software \cite{elman2014ifiss} for the steady-state solution of a linearized 2D Navier-Stokes equations and are obtained from  for the lid-driven cavity problem with a spatial discretization based on $Q_2/Q_1$ elements and a uniform grid, varying the Reynolds number between $100$ and $900$.  Since our method is a black-box algebraic approach, we ignore the fact that these problems arise from specific applications and treat them just as given test matrices. It is clear that problem specific solvers, see \cite{ElmSW14}, that exploit the underlying discrete differential properties of the operators will, in general, show better performance than black-box methods.  The $B$ matrix in this case  is numerically rank-deficient. Therefore,  in the new scheme, we compute a nullspace basis of size $m-n+1$. In Table~\ref{tab:generalized:matrices}, we provide the smallest and largest eigenvalues of $A+A^T$. For three of the test problems in this set, $A+A^T$ is positive definite but close to singular, while for the remaining five problems it is robustly positive definite.

For the general case, nine test matrices are used (given in Table~\ref{tab:general:matrices}). Six matrices are obtained from the SuiteSparse Matrix Collection and three are randomly generated. Here, garon1 and garon2  arise in the 2D finite element discretization of the  Navier-Stokes  equations with inlet and outlet on opposing sides of the geometry. Tolosa matrices (Tols90-4000)  arise in the stability analysis of an airplane flight model.  The three random matrices are generated as follows: Each block is generated randomly, and to avoid a singular system, a small perturbation of the form $\xi I_n$ is added to the $A$ block and $\xi [I_m\ 0]^T$ to the blocks $C$ and $B$, i.e. for given $m$ and $n$,   $A\leftarrow \xi I_{n} +  R_{n\times n}$,  $B\leftarrow \xi[I_{m}\ 0]^T +  R_{n \times m} $ and
$C^T\leftarrow \xi[I_{m}\ 0]+  R_{m \times n} $ where $R_{n\times m}$ is a sparse random $n\times m$ matrix with uniformly distributed nonzeros between $0$ and $1$,  with a density of $1\%$. In our tests, we chose $\xi=10^{-1}$.  For all cases, if the right-hand side vector is available, we used it; otherwise, we generated the right-hand side vector from a solution vector of all ones to have a consistent system.
\begin{table}[htbp]
\begin{tabular}{ll|rrrll}
\#& Matrix Name &n & m & nnz & $\lambda_{min}(A)$ & $\lambda_{max}(A)$\\\hline
1 & reorientation\_1 & 396 & 281 &7,326 & $-1.7\times10^6$ & $1.0\times 10^9$\\
2 & reorientation\_4 &1,596 & 1,121 &33,630 & -43.2 & $3.8\times10^8$ \\
3 & reorientation\_8 &1,826 & 1,282 &37,894 & -43.3  & $5.9\times10^8$ \\
4 & dynamicSoaringProblem\_1 &363 & 284  &5,367 & $7.1\times10^{-7}$ & 237.1 \\
5 & dynamicSoaringProblem\_4 &1,794 & 1,397 & 36,516 & $1.1\times10^{-8}$ & $4.2\times 10^5$ \\
6 & goddardRocketProblem\_2 &434 & 433 & 9,058 & $-1.1\times 10^3$ & $1.1\times 10^3$ \\
7 & tumorAntiAngiogenesis\_1 &123 & 82 &2,748 & -25.3 &  65.2\\
8 & tumorAntiAngiogenesis\_8 & 294 & 196 &4,776 & -99.6 & $7.3\times 10^{11}$
\end{tabular}
\caption{Symmetric case: test problems}
\label{tab:symmetric:matrices}
\end{table}

\begin{table}[htbp]
\begin{tabular}{ll|rrrll}
\#& Matrix Name &n & m & nnz & $\lambda_{min}(A+A^T)$ & $\lambda_{max}(A+A^T)$\\\hline
1 & rajat04 & 1,008 & 33 &8,725 & 1.1 & $2.1\times 10^4 $ \\
2 & rajat14 & 171 & 9 &1,475 & $5.8\times10^{-3}$ & $8.6\times10^4$\\
3 & fpga\_trans\_01 &1,154 & 66 & 7,382 & $8.7\times 10^{-10}$ & 2.6\\
4 & drivencavity\_Re100 & 578 & 81  &10,814 & $7.6\times 10^{-4}$ & 1.0 \\
5 & drivencavity\_Re200 & 578 & 81  &10,814 & $3.3\times10^{-4}$ & 1.0 \\
6 & drivencavity\_Re500 & 578 & 81  &10,814  & $-8.2\times 10^{-4}$ & 1.0\\
7 & drivencavity\_Re700 & 578 & 81  &10,814 & $-2.1\times 10^{-3}$ & 1.0 \\
8 & drivencavity\_Re900 & 578 & 81  &10,814 & $-9.1\times 10^{-3}$ & 1.0 \\
\end{tabular}
\caption{Structurally symmetric case: matrix properties}
\label{tab:generalized:matrices}
\end{table}

\begin{table}[htbp]
\begin{tabular}{ll|rrrr}
\#& Matrix Name &n & m & nnz \\\hline
1 & garon1 & 2,775 & 400 &84,723 \\
2 & garon2 & 11,935 & 1,600 &373,235 \\
3 & tols90 &72 & 18 & 1,746\\
4 & tols340 & 272 & 68  &2,196 \\
5 & tols1090 & 872 & 218  &3,546 \\
6 & tols4000 & 3,200 & 800  &8,784 \\
7 & random1 (cond.$2.15\times 10^4)$ & 100 & 90  &2,900\\
8 &random2 (cond.$2.89\times 10^9$ )& 100 & 90  & 554\\
9 &random3 (cond.$2.26\times 10^4$ )& 1,000 & 900  & 3,079\\
\end{tabular}
\caption{General case: matrix properties}
\label{tab:general:matrices}
\end{table}

\section{Experimental framework}
\label{sec:framework}
The following numerical results were obtained using MATLAB R2023b on an Apple Macbook Air with an M2 processor and 8 GB of memory. As a baseline of comparison, we use restarted GMRES(m) with $m=10$, preconditioned with MATLAB's incomplete $LU$ method with pivoting and thresholding, ILUTP($\tau$), where $\tau$ is the drop tolerance.  We apply reverse Cuthill-McKee reordering to the systems before computing the incomplete factors by ILUTP($\tau$)  to improve the robustness of the baseline method, see \cite{benzi1999orderings}.  Here, GMRES(m)  is preconditioned from the left and right.
We note that for the symmetric case, an incomplete $LDL^T$ factorization based preconditioner would be more effective and suitable as a baseline, see~\cite{greif2017sym}.
 Therefore, we use ILUTP($\tau$) for all cases, including the symmetric case. Then in the symmetric case we also use MINRES, preconditioned with the \emph{absolute value preconditioner}, $L|D|L^T$,  obtained by modifying the  incomplete Bunch-Kaufman-Parlett factorization \cite{greif2017sym}, \emph{sym-ildl}, where,  as in \cite{vecharynski2013absolute},  $|D|=V|\Lambda|V^T$ if $D$ has a spectral decomposition $V\Lambda V^T$.

For our new schemes, we use the restarted and right preconditioned flexible GMRES; fGMRES(m), with $m=10$, as the solver for the outer iteration. We chose fGMRES(m) because the proposed preconditioners (comprised of other iterative solvers)  inherently change between iterations. The new schemes have a number of parameters: (a) drop tolerances and thresholds: $\tau_{SAROC}, \rho_{SAROC}, \tau_{FSAI}, \rho_{FSAI}$
and (b) Stopping tolerances for the outer fGMRES, inner LSQR, inner fGMRES and innermost MRS (or inner GC). Furthermore, if $M$-orthogonalization is performed on the nullspace basis, we include the window size ($w_{MGS}$) and drop tolerance ($\tau_{MGS}$) parameters. To illustrate the robustness of the new solver class, we experiment with small, large, and mixed thresholds and tolerances. These are  given in Table~\ref{tab:parameters}.  For ILUTP($\tau$),  $\tau$ is set to  $10^{-4},10^{-3},10^{-2}$ and $10^{-1}$.  {For \emph{sym-ildl}, to compute the incomplete factors,  we use the default parameters. That is, the \emph{drop tolerance} and \emph{fill-in} are  $10^{-3}$ relative to the norm of the column  and at most $3$ times the average number of nonzeros per row, respectively. Furthermore, \emph{Rook} pivoting and  \emph{AMD} reordering are used. We refer the reader to \cite{greif2017sym} for details of the implementation of these parameters.} One of the main difficulties for classical incomplete factorization-based preconditioners is that there is no easy way to predict a good drop tolerance that would work well for a given problem. Sometimes, it is even counterintuitive that decreasing the drop tolerance might result in encountering a zero pivot.

For  the new schemes and the baseline methods,  we start the iterations with the zero vector and stop when the $L_2$ norm of the  relative residual drops below  $10^{-5}$ or the maximum number (1,000) of iterations is reached.
Note that for the case that good starting vectors are available (for example if one has a sequence of slightly changing matrices or right hand sides) then the convergence will be substantially faster.

The residual norms are obtained as a by-product of fGMRES and GMRES.  {Here,} right preconditioned fGMRES produces the \emph{true residual} norm, and GMRES produces the \emph{preconditioned residual} norm due to ILUTP-based left and right preconditioning. For these reasons, we also calculate and report the final \emph{true relative residual} norms.  For the new schemes,  { for the symmetric case, the inner LSQR and CG iterations are stopped when the $L_2$-norm of the  relative residual is less than or equal to $\epsilon_{in}$ and $\epsilon_{innermost}$  is not used. } For the structurally symmetric and general case, inner LSQR and fGMRES iterations are stopped when the $L_2$-norm of the relative residual is less than or equal to $\epsilon_{in}$ and the innermost MRS iterations are stopped when the $L_2$-norm of the relative residual is less than or equal to $\epsilon_{innermost}$. Since our new schemes are based on multiple nested loops and external functions,  a fair comparison in terms of the running times against MATLAB's built-in functions is currently not possible, and we only provide comparisons in terms of the number of iterations, number of nonzeros, and the robustness of the new and baseline schemes. The number of iterations and the number of nonzeros provide approximate measures for computational time and memory requirements.

\begin{table}[htbp]\centering
\begin{tabular}{r|lll}
&\multicolumn{3}{c}{Tolerances}  \\
 & $\tau_{large}$ & $\tau_{mix}$ & $\tau_{small}$ \\ \hline
$\tau_{SAROC}$ & $10^{-3}$ & $10^{-2}$ & $10^{-5}$ \\
$\rho_{SAROC}$ & $10^{-3}$ & $10^{-2}$ & $10^{-5}$ \\
$\tau_{FSAI}$ & $10^{-3}$ & $10^{-3}$ & $10^{-5}$ \\
$\rho_{FSAI}$ & $10^{-3}$ & $10^{-3}$ & $10^{-5}$\\ \hdashline
$\tau_{MGS}$ & $10^{-3}$ & $10^{-2}$ & $10^{-5}$ \\
$w_{MGS}$ & $5$ & $5$ & $15$ \\\hdashline
$\epsilon_{in}$ & $10^{-3}$ & $10^{-4}$ &  $10^{-5}$\\
$\epsilon_{innermost}$ & $10^{-3}$ & $10^{-5}$ &  $10^{-5}$\\
\end{tabular}
\caption{Parameters of the new schemes}
\label{tab:parameters}
\end{table}

\section{Numerical results}
\label{sec:results}
In this section we present the numerical results obtained by our new solver class. Throughout the presented results, we denote by $\star$  if the true residual is  greater than the stopping criterion, by $\ddagger$ if the maximum number of iterations is reached, and  by $\S$ if MATLAB runs out of memory. Furthermore, for the baseline (ILUTP) scheme, $\dagger$ indicates a failure in the incomplete factorization, since a zero pivot is encountered. All average values are rounded to one digit after the decimal point.

 { \subsection{Symmetric case}
First, we consider the new scheme without $M$-orthogonalization.  For the problems given in Table~\ref{tab:symmetric:matrices}, we present in Table~\ref{tab:symmeric:results:outer} the number of outer GMRES(10), MINRES and fGMRES(10) iterations for the baselines and the new scheme, respectively. As expected, the smaller the drop tolerances, the more likely the new method is to converge. The same is true for ILUTP, however, for three cases (reorientation\_4, reorientation\_8 and tumorAntiAngiogenesis\_8), the final true relative residual norms are below the stopping tolerance even though the preconditioned relative residuals norms indicate that the method has converged. The encounter of a zero pivot is a frequent mode of failure for ILUTP, while the new scheme only failed by reaching the maximum allowed number of iterations without converging to the required stopping tolerance. In Table~\ref{tab:symmetric:results:memory}, we report the total number of nonzeros for the new scheme (i.e., $nnz(\widetilde{Z})+nnz(\widetilde{W})$), ILUTP($\tau$) (i.e., $nnz(L)+nnz(U)$)  and sym-ildl (i.e., $nnz(L)+nnz(D)$). The best method among those that converge without any failure for each problem is given in {\bf bold}. In all tests  and in all three cases, the new scheme  achieved fewer nonzeros compared to ILUTP and sym-ildl, respectively.

In Table~\ref{tab:symmeric:results:inner}, we present the average number of CG and LSQR iterations for the new scheme. Although the number of CG iterations is reasonably small for all cases, the number of LSQR iterations is relatively large, especially for small drop and threshold tolerances. This is to be expected, since we do not use any further preconditioning for the least squares problems, while CG is already preconditioned with SPAI. We note that for two cases (dynamicSoaringProblem\_4 and godddardRocketProblem\_2), even though the LSQR fails to converge within $1,000$ iterations, the outer iterations still converge without any failure.
\begin{table}[htbp]\centering
\begin{tabular}{r|lll|llll|l}
&\multicolumn{3}{c}{New Scheme} &\multicolumn{4}{c}{ILUTP($\tau$)} &  \\
\# & $\tau_{large}$ & $\tau_{mix}$ & $\tau_{small}$  & $10^{-4}$ & $10^{-3}$ &$10^{-2}$  & $10^{-1}$ & sym-ildl\\ \hline
1 & 17 & $\ddagger$ & 2 & $\dagger$ &$\dagger$ &$\dagger$ &$\dagger$ & $\ddagger$ \\
2 & $\ddagger$ & $\ddagger$ & 2 & $\dagger$ &$\dagger$ &$\dagger$ &$\dagger$ & 31\\
3 & 718 & $\ddagger$ & 2 & $\dagger$  &$\dagger$ &$\dagger$ &$\dagger$ & 29 \\
4 & 14 & 23 & 2 & 3 & 5  &$\dagger$ &$\dagger$ & 83\\
5 & 99 & $\ddagger$ &  4  & $6^{\star}$ & $30^{\star}$  &$\dagger$ &$\dagger$ & 117 \\
6 & $\ddagger$ & $\ddagger$ &  12  & 7 & $\ddagger$  &$\dagger$ &$\dagger$ & $\ddagger$\\
7 & 3 & 6 &  2  & 4 & 6  & 9 &$\dagger$ & 125 \\
8 & $\ddagger$ & $\ddagger$ & 1 & $3^{\star}$ & $7^{\star}$ &$\dagger$ &$\dagger$ & 3 \\
\end{tabular}
\caption{Number of outer fGMRES/GMRES/MINRES iterations for the symmetric case}
\label{tab:symmeric:results:outer}
\end{table}
\begin{sidewaystable}[htbp]\centering
\begin{tabular}{r|lll|llll|l}
&\multicolumn{3}{c}{New Scheme} &\multicolumn{4}{c}{ILUTP($\tau$)} \\
\# & $\tau_{large}$ & $\tau_{mix}$ & $\tau_{small}$  &$10^{-4}$ & $10^{-3}$ &$10^{-2}$  & $10^{-1}$ & sym-ildl  \\ \hline
1 & {\bf 21,512} & 1,065$^\ddagger$ & 37,526 &$\dagger$ &$\dagger$ &$\dagger$ &$\dagger$ & 3,639$^\ddagger$ \\
2 & 81,544$^\ddagger$ & 1,073$^\ddagger$ & 572,767& $\dagger$ &$\dagger$ &$\dagger$ &$\dagger$ & {\bf 52,583} \\
3 & 93,634 & 1,353$^\ddagger$ & 686,656 & $\dagger$ &$\dagger$ &$\dagger$ &$\dagger$ &  {\bf 58,407}  \\
4 & 23,660 & 15,623  & 25,581& 64,783 & 49,684  &$\dagger$ &$\dagger$ & {\bf 9,626} \\
5 & 480,578 & 61,137$^\ddagger$ &  631,284 & 750,882$^{\star}$ & 330,140$^{\star}$  &$\dagger$ &$\dagger$ & {\bf 55,834} \\
6 & $92^\ddagger$ & $44^\ddagger$ &  {\bf 93} & 96,486 & 67,828$^{\ddagger}$  &$\dagger$ &$\dagger$ & 14,535$^\ddagger$\\
7 & 2,698 & {\bf 1,502} &  4,203 & 17,769 & 14,439  & 8,017  & $\dagger$ & 3,131\\
8 & 2,987$^\ddagger$ & 1,318$\ddagger$ & 9,217 & 57,502$^{\star}$ & 39,290$^{\star}$ &$\dagger$ &$\dagger$ & {\bf 6,775} \\
\end{tabular}
\caption{Number of nonzeros of the preconditioners for the symmetric case}
\label{tab:symmetric:results:memory}
\end{sidewaystable}
\begin{table}[htbp]\centering
\begin{tabular}{r|ll|ll|ll}
\# &\multicolumn{2}{c}{$\tau_{large}$} &\multicolumn{2}{c}{$\tau_{mix}$}  & \multicolumn{2}{c}{$\tau_{small}$}  \\
 & cg & lsqr & cg & lsqr & cg & lsqr \\ \hline
1 & $2.9$ & $238.5$ & $2.0$ & $646.8$ & $1.5$ & $874.8$ \\
2 & $8.1$ & $145.9$ & $1.9$ & $982.1$ & $9.0$ & $767.5$ \\
3 & $5.8$ & $153.7$ & $1.9$ & $992.4$ & $20.5$ & $766.3$ \\
4 & $2.0$ & $84.2$ & $3.0$ & $350.8$ & $1.0$ & $602.5$ \\
5 & $2.2$ & $272.9$ & $4.3$ & $881.9$ & $2.0$ & $\ddagger$ \\
6 & $1.0$ & $71.4$ & $1.0$ & $520.1$ & $1.0$ & $\ddagger$ \\
7 & $1.0$ & $138.7$ & $2.0$ & $140.6$ & $1.0$ & $144.5$ \\
8 & $1.4$ & $115.0$ & $6.0$ & $518.1$ & $5.0$ & $768.5$ \\
  \end{tabular}
  \caption{Average number of inner CG and LSQR iterations for the symmetric case}
\label{tab:symmeric:results:inner}
\end{table}
Figure~\ref{fig:truerelativeresidual:symmetric} shows the relative residual norms of the new scheme and the baseline  for the symmetric case. The  horizontal dashed line visualizes the stopping tolerance where the norm of the relative residual is  equal to $10^{-5}$.
\begin{figure}[htbp]
\centering
\includegraphics[width=1\linewidth]{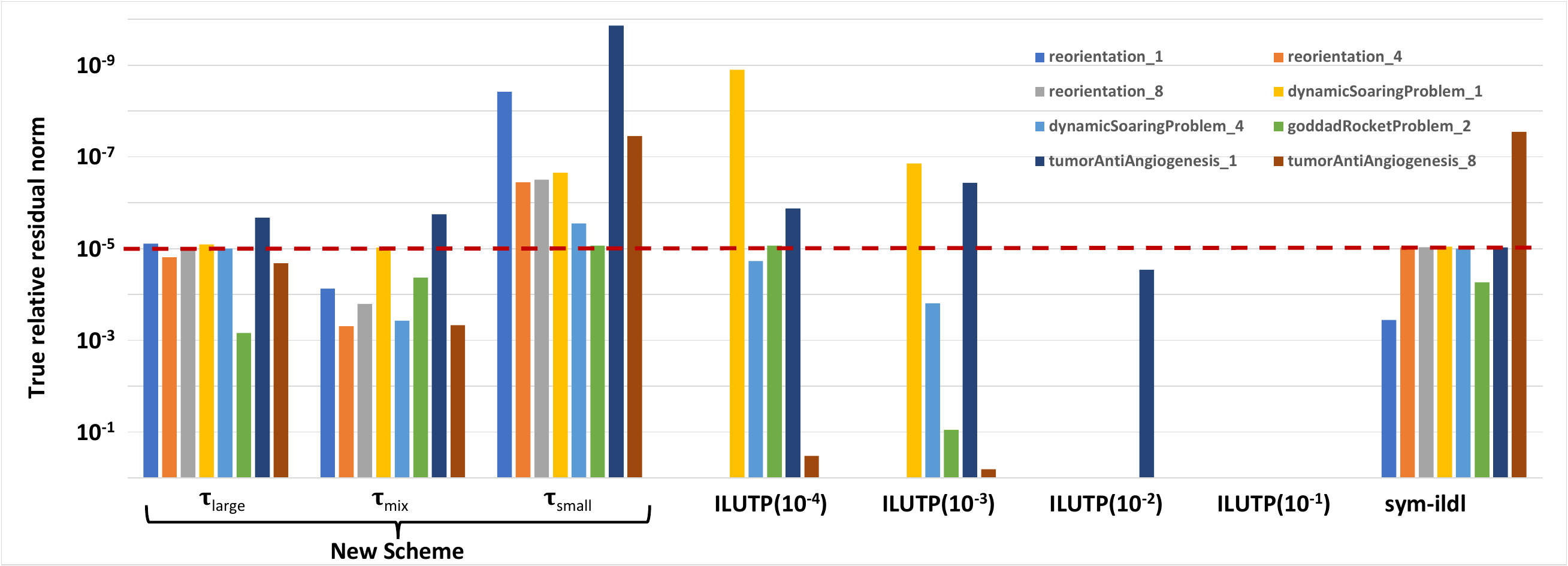}
  \caption{Final true relative residual norms of the new scheme and the baseline method for the symmetric case}
\label{fig:truerelativeresidual:symmetric}
\end{figure}

Next, we consider the additional $M$-orthogonalization of the nullspace basis via the modified Gram-Schmidt procedure. The number of outer-inner iterations and the total number of nonzeros for the new scheme are presented in Tables~\ref{tab:symmetric:results:outer:morth}, \ref{tab:symmeric:results:inner:morth} and \ref{tab:symmetric:results:nnz:morth}, respectively. In this case, the total number of nonzeros refers to the sum of the number of  nonzeros of the nullspace basis, the $M$-orthogonalized nullspace basis  and  the nonzeros of the SPAI preconditioner (i.e., $nnz(\widetilde{Z})+nnz(\bar{Z})+nnz(\tilde{W})$). Although $\widetilde{Z}$ can be discarded after $\bar{Z}$ is computed, we include it in the sum as a measure of the peak memory requirement.  When $\tau_{small}$ is used, there is not much difference in terms of inner and outer iterations, except for reorientation\_8, where the number of average CG iterations drops from $20.5$ to $2.5$ if the nullspace basis is $M$-orthogonalized. When $\tau_{mix}$ is used, there are two problems for which the new scheme did not converge without $M$-orthogonalization, namely  goddardRocketProblem\_2 and tumorAntiAngiogenesis\_8, these converge now with $M$-orthogonalization.  GoddardRocketProblem\_2, is particularly interesting, since it has a nullspace dimension of $1$, and hence $M$-orthogonalization is just normalizing the basis vector. 
The matrix $A$ has a condition number of $1.85\times 10^{15}$, i.e. it is almost singular, even though the condition number of the coefficient matrix $\mathcal W$ is $1.38\times 10^{7}$. In this case, the nullspace method becomes a more effective preconditioner with $M$-orthogonalization and the number of outer fGMRES iterations is reduced to $58$, without any increase in memory requirement.    When $\tau_{large}$ is used, $M$-orthogonalization appears to cause more failures, because the nullspace basis is likely to be inaccurate to begin with, and further approximations only increase the error.
\begin{table}[htbp]\centering
\begin{tabular}{r|lll}
\# & $\tau_{large}$ &$\tau_{mix}$ & $\tau_{small}$ \\\hline
1 & 19 & $\ddagger$   &  2 \\
2 & $\ddagger$  & $\ddagger$   & $\S$  \\
3 & $\ddagger$ &  $\ddagger$ & 2   \\
4 & 12 &  20   & 2    \\
5 & 570 & $\ddagger$  & 4   \\
6 & $\ddagger$ &  58 &  12  \\
7 & 3  &  6  &  2  \\
8 & $\ddagger$ & 1  & 1
  \end{tabular}
  \caption{Number of outer fGMRES iterations for the new scheme (symmetric case) with $M$-orthogonalization}
\label{tab:symmetric:results:outer:morth}
\end{table}
\begin{table}[htbp]\centering
\begin{tabular}{r|ll|ll|ll}
\# &\multicolumn{2}{c}{$\tau_{large}$} &\multicolumn{2}{c}{$\tau_{mix}$}  & \multicolumn{2}{c}{$\tau_{small}$}  \\
 & cg & lsqr & cg & lsqr & cg & lsqr \\ \hline
1 & 2.0 & 235.9 &2  & 641.8 & 1.0 & 866.5  \\
2 & 2.0 & 157.3 & 1.9 & 980.6  & $\S$ & $\S$ \\
3 & 2.4 & 153.0 & 1.9 & 958.3 & 2.5 & 766.3 \\
4 & 2.0 & 87.1 & 2.0 & 354.0 &1.0  & 599.5 \\
5 & 2.0 & 285.5 &3.0  & 848.9 & 1.8 & $\ddagger$ \\
6 & 1.0 & 68.5 & 1.0 & 511.8 &1.0  & $\ddagger$ \\
7 & 2.0 & 136.0 & 2.0 &141.2  & 1.0 & 144.8 \\
8 & 8.0 & 128.6 & 4.0 & 721.0 & 6.0 & 763.0 \\
  \end{tabular}
  \caption{Average number of inner CG and LSQR iterations  for the new scheme (symmetric case) with $M$-orthogonalization}
\label{tab:symmeric:results:inner:morth}
\end{table}
\begin{table}[htbp]\centering
\begin{tabular}{r|lll}
\# & $\tau_{large}$ &$\tau_{mix}$ & $\tau_{small}$ \\\hline
1 & 36,413 & 1,642  & 73,858 \\
2 & 142,321 &1,590 & $\S$ \\
3 & 162,901 &2,061 & 1,272,471 \\
4 & 45,250 & 27,759& 50,769 \\
5 & 881,668 & 94,735 & 1,218,169 \\
6 & 183 & 87 & 185 \\
7 &4,681  &  2,276 & 8,145  \\
8 &6,045  & 1,952& 24,301 \\
  \end{tabular}
  \caption{Number of nonzeros of the new scheme (symmetric case) with $M$-orthogonalization}
\label{tab:symmetric:results:nnz:morth}
\end{table}
}

\subsection{Structurally symmetric case}
As before, let us  first consider the new scheme without $M$-ortho\-go\-na\-li\-za\-tion.  For the problems given in Table~\ref{tab:generalized:matrices}, in Table~\ref{tab:generalized:results:outer} we present the number of outer GMRES(10) and fGMRES(10) iterations for the baseline method and the new scheme, respectively. Here, ILUTP frequently fails, as it encounters a zero pivot, sometimes even when a stringent or a relaxed drop tolerance is used but converges for some arbitrary values of drop tolerances for the lid-driven cavity problem. In a number of cases, ILUTP fails because the final true relative residual norm is larger than the stopping tolerance. The new method converges for all problems without any failure, showing the robustness of the new scheme for the structurally symmetric problem set.
\begin{table}[htbp]\centering
\begin{tabular}{r|lll|llll}
&\multicolumn{3}{c}{New Scheme} &\multicolumn{4}{c}{ILUTP($\tau$)} \\
\# & $\tau_{large}$ & $\tau_{mix}$ & $\tau_{small}$  & $10^{-4}$ & $10^{-3}$ &$10^{-2}$  & $10^{-1}$ \\ \hline
1 & 3 & 2 & 2 &$\dagger$ &$\dagger$ &$\dagger$ &$\dagger$ \\
2 & 2 & 2 & 2 & $\dagger$ &$\dagger$ &$\dagger$ &$\dagger$ \\
3 & 2 & 1 & 1  & 4 & 6  &12  &13  \\
4 & 2 & 2 & 2 &$2^\star$ & 3  & 6 & $32^\star$ \\
5 & 2 & 3 &  2 &$1^\star$ & 3  &7 &$44^\star$ \\
6 & 3 & 3 &  1  & 2 & 3  & 7 &$\ddagger$ \\
7 & 3 & 3 &  2  & $\dagger$ & $\ddagger$  & $9^\star$ &$\ddagger$ \\
8 & 4 & 3 & 2  & $\dagger$ & $\ddagger$  & 10 &$\dagger$ \\
\end{tabular}
\caption{Number of outer fGMRES/GMRES iterations for the structurally symmetric case}
\label{tab:generalized:results:outer}
\end{table}

\begin{figure}[htbp]
\centering
  \includegraphics[width=1\linewidth]{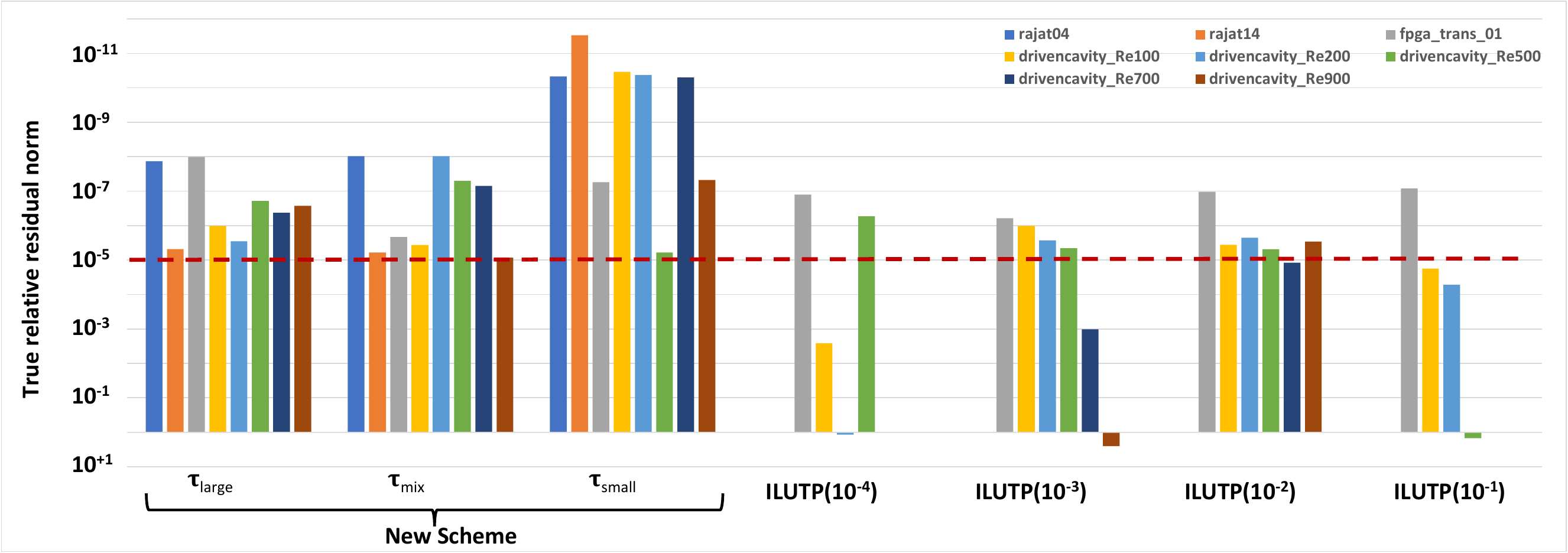}
  \caption{Final true relative residual norms of the new scheme and the baseline for the structurally symmetric case}
  \label{fig:truerelativeresidual:generalized}
\end{figure}

In Table~\ref{tab:generalized:results:memory}, we report the total number of nonzeros for the new scheme (i.e., $nnz(\widetilde{Z})+nnz(\widetilde{W})$)  and ILUTP($\tau$) (i.e., $nnz(L)+nnz(U)$). The best method among those that converge without any failure for each problem is given in {\bf bold}. The new scheme achieves fewer nonzeros compared to ILUTP for 3 cases, while ILUTP is better than the new scheme for 5 cases. However, in most cases, the difference in the number of nonzeros is not significant.
\begin{sidewaystable}[htbp]\centering
\begin{tabular}{r|lll|llll}
&\multicolumn{3}{c}{New Scheme} &\multicolumn{4}{c}{ILUTP($\tau$)} \\
\# & $\tau_{large}$ & $\tau_{mix}$ & $\tau_{small}$ & $10^{-4}$ & $10^{-3}$ &$10^{-2}$  & $10^{-1}$ \\ \hline
1 & {\bf 23,856} & {\bf 23,856} & 203,761 &$\dagger$ &$\dagger$ &$\dagger$ &$\dagger$ \\
2 & {\bf 1,329} & {\bf 1,329} & 2,326 & $\dagger$ &$\dagger$ &$\dagger$ &$\dagger$ \\
3 & 10,946 & 10,946 & 26,596 & 11,526 & 8,521 &4,905  &{\bf 4,417} \\
4 & 55,661 & 51,984 & 69,923 &66,804$^\star$ & 58,686  & {\bf 41,156} & 16,444$^\star$ \\
5 & 55,591 & 55,584 &  69,946 &68,323$^\star$ & 60,242  &{\bf 46,078} & 19,095$^\star$\\
6 & 58,271 & 58,266 &  70,325 & 67,989 & 60,271  & {\bf 45,570} &20,011$^\ddagger$ \\
7 & {\bf 60,019} & 60,042 & 70,842  & $\dagger$ & 64,229$^\ddagger$  & 46,681$^\star$ &  $\dagger$ \\
8 & 63,143 & 63,118 & 71,699  & $\dagger$ & 62,842$^\ddagger$  & {\bf 45,280} &  $\dagger$ \\
\end{tabular}
\caption{Number of nonzeros of the preconditioners for the structurally symmetric case}
\label{tab:generalized:results:memory}
\end{sidewaystable}
In Table~\ref{tab:generalized:results:inner}, we present the average number of inner fGMRES(10),  LSQR and MRS iterations for the new scheme. While all are reasonably small for the circuit simulation systems, the number of inner fGMRES(10) and MRS iterations grow as the Reynold's number increases for the lid-driven cavity problem. However, this growth is suppressed if $\tau_{large}$ is used. For structurally symmetric test problems, the number of LSQR iterations remained reasonably small for all problems, even without using any preconditioner.  Figure~\ref{fig:truerelativeresidual:generalized} shows the final true relative residual norms of the new scheme and the baseline for the structurally symmetric case. The  horizontal dashed line visualizes the stopping tolerance where the norm of the relative residual is equal to $10^{-5}$.
\begin{table}[htbp]\centering
\begin{tabular}{r|lll|lll|lll}
\# &\multicolumn{3}{c}{$\tau_{large}$} &\multicolumn{3}{c}{$\tau_{mix}$}  & \multicolumn{3}{c}{$\tau_{small}$}  \\
 & fgmres & lsqr & mrs & fgmres  & lsqr & mrs & fgmres  & lsqr & mrs \\ \hline
1 &3.0 & 1.5  & 2.9 & 4.0 & 1.5  & 4.0 & 2.0 & 1.5 & 4.0  \\
2 &2.0 & 1.0 & 5.0 & 2.0 & 1.0 &  7.0 & 2.0& 1.0 & 7.0\\
3 & 3.0 & 1.0 & 2.0 & 4.0 & 1.0 & 3.3 & 2.0 & 1.0 & 3.5\\
4 & 2.5 & 16.8 & 17.8 & 3.0 & 24.3 & 27.5 & 2.0& 25.3&27.5\\
5 & 2.0 & 16.5 & 31.8 & 3.0 & 23.7 & 48.0 & 2.0 & 24.8&  47.8\\
6 & 3.7 & 17.3 & 101.4 & 5.0 & 23.3 & 164.5& 3.0 & 23.5 & 151.7\\
7 & 9.0 & 17.5 & 230.3 & 10.0 & 23.8 & 374.0&7.0 & 23.8 & 360.7\\
8 &52.5 & 18.3 & 493.9 & 81.3 & 22.8 & 791.9& 132.5& 27.8 & 809.6
  \end{tabular}
  \caption{Average number of inner fGMRES,  LSQR and MRS iterations for the new scheme (structurally symmetric case)}
\label{tab:generalized:results:inner}
\end{table}

Next, we consider the additional $M$-orthogonalization of the nullspace basis via the modified Gram-Schmidt procedure. The number of outer fGMRES(10) iterations and the total number of nonzeros for the new scheme are given in Tables~\ref{tab:generalized:results:outer:morth} and \ref{tab:generalized:results:nnz:morth}, respectively. As before, the total number of nonzeros refers to the sum of the number of nonzeros of the nullspace basis, the $M$-orthogonalized nullspace basis and the nonzeros of the SPAI preconditioner (i.e. $nnz(\widetilde{Z})+nnz(\bar{Z})+nnz(\tilde{W})$).  The average number of inner fGMRES(10), LSQR, and MRS are given in Table~\ref{tab:generalized:results:inner:morth}. For test problems in the structurally symmetric case, we do not observe any significant differences when $M$-orthogonalization is performed. There is some increase in the number of nonzeros in most cases.
\begin{table}[htbp]\centering
\begin{tabular}{r|lll}
\# & $\tau_{large}$ &$\tau_{mix}$ & $\tau_{small}$ \\\hline
1 & 2 &  2  &  1 \\
2 & 2  & 2   & 2  \\
3 & 2 & 1  & 1   \\
4 &3  &  4   & 2    \\
5 & 3 & 4  & 2   \\
6 & 3 &  3 & 2   \\
7 & 3  & 5   & 2   \\
8 & 4 &  5 & 2
  \end{tabular}
  \caption{Number of outer fGMRES iterations for the new scheme (structurally symmetric case) with $M$-orthogonalization}
\label{tab:generalized:results:outer:morth}
\end{table}
\begin{table}[htbp]\centering
\begin{tabular}{r|lll}
\# & $\tau_{large}$ &$\tau_{mix}$ & $\tau_{small}$ \\\hline
1 & 25,850 & 25,640  & 194,360 \\
2 & 1,510 & 1,458  & 2,522 \\
3 &12,408 & 12,452  & 27,714   \\
4 & 62,499 &  61,007 & 94,424    \\
5 & 62,508 & 60,989  &  94,444   \\
6 & 64,497 &  63,109 &  95,064  \\
7 & 65,518 & 66,266  & 95,916   \\
8 & 70,609 & 77,180  & 98,515
  \end{tabular}
  \caption{Number of nonzeros of the new scheme for the new scheme (structurally symmetric case) with $M$-orthogonalization}
\label{tab:generalized:results:nnz:morth}
\end{table}
\begin{table}[htbp]\centering
\begin{tabular}{r|lll|lll|lll}
\# &\multicolumn{3}{c}{$\tau_{large}$} &\multicolumn{3}{c}{$\tau_{mix}$}  & \multicolumn{3}{c}{$\tau_{small}$}  \\
 & fgmres & lsqr & mrs & fgmres  & lsqr & mrs & fgmres  & lsqr & mrs \\ \hline
1 & 3.0 & 1.5   & 2.8  & 3.5 & 1.5  & 4.0 & 2.0 & 1.5 &  4.0  \\
2 & 2.0 &  1.0  &5.0  &  2.0 & 1.0  & 7.0 & 1.5 & 1.0 & 7.0  \\
3 & 2.5 & 1.0  &  2.0 & 3.0 & 1.0  & 3.3 & 2.0 & 1.0 & 3.0  \\
4 &2.0 & 19.8   & 18.5 & 3.0 & 22.1  & 27.8 & 2.0 & 29.5 & 27.5  \\
5 & 2.3 & 20.0  & 32.0 & 3.0 & 21.8  & 48.4 & 2.0 & 29.0 & 48.5  \\
6 & 3.0 & 18.2  &100.1  & 4.8 & 21.6  & 164.8 & 3.0 & 29.0  & 151.7  \\
7 & 8.0 & 19.2  & 216.0 & 54.8 & 21.4  & 567.2 & 7.0 & 29.3 & 358.4  \\
8 & 54.5 & 17.9  & 460.9  & 72.0 & 21.4  & 752.2  & 124.5 & 28.0 & 810.6
  \end{tabular}
  \caption{Average number of inner fGMRES,  LSQR and MRS iterations for the new scheme  (structurally symmetric case) with $M$-orthogonalization}
\label{tab:generalized:results:inner:morth}
\end{table}

\subsection{General case}
 For the general problems given in Table~\ref{tab:general:matrices},  in  Table~\ref{tab:general:results:outer} we present the number of outer GMRES(10) and fGMRES(10) iterations for the baseline method and the new scheme, respectively. The new scheme converges for all cases where $\tau_{mix}$ and $\tau_{small}$ are used, while ILUTP has at least one failure due to the encounter of a zero pivot for all drop tolerances. The new scheme fails in one case when $\tau_{large}$ is used due to reaching the maximum allowed iterations without reaching the required relative residual norm.
\begin{table}[htbp]\centering
\begin{tabular}{r|lll|llll}
&\multicolumn{3}{c}{New Scheme} &\multicolumn{4}{c}{ILUTP($\tau$)} \\
\# & $\tau_{large}$ & $\tau_{mix}$ & $\tau_{small}$ &  $10^{-4}$ & $10^{-3}$ &$10^{-2}$  & $10^{-1}$ \\ \hline
1 & 3 & 3 & 2 & 5 &$\dagger$ &$\dagger$ &$\dagger$ \\
2 & 10 & 9 & 2  & 9 &$\dagger$ &$\dagger$ &$\dagger$ \\
3 & 2 & 2 & 1  & 2 & 2  & 3  & 4  \\
4 & 7 & 2 & 1 & 2 & 2  & 3 &  4 \\
5 & 34 & 5 &  4  & 2 & 2  & 3 &4  \\
6 & 27 & 11 & 11  & 2 & 2  & 3 & 4 \\
7 & 4 & 17 &  1 & 2 & 4  & 10 &$\dagger$ \\
8 & $\ddagger$  & 53 & 5 & $\dagger$ & $\dagger$  & $\dagger$ &$\dagger$ \\
9 & 9 & 2 & 2  & 1 &  $\dagger$  & $\dagger$ &$\dagger$ \\
\end{tabular}
\caption{Number of outer fGMRES/GMRES iterations for the general case}
\label{tab:general:results:outer}
\end{table}
In Table~\ref{tab:general:results:memory}, we report the total number of nonzeros for the new scheme (i.e., $nnz(\widetilde{Z})+nnz(\widetilde{U})+nnz(\widetilde{W})$)  and ILUTP($\tau$) (i.e., $nnz(L)+nnz(U)$). The best method among those that converge without any failure for each problem is given in {\bf bold}.  The new scheme not only converges without any failures for all problems but also achieves fewer nonzeros compared to ILUTP for all cases, except for two (garon\_1 and garon\_2). We note that for these two cases, even though ILUTP achieves fewer number of nonzeros, it is likely to fail if the drop tolerance is not chosen correctly.

\begin{figure}[htbp]
\centering
  \includegraphics[width=1\linewidth]{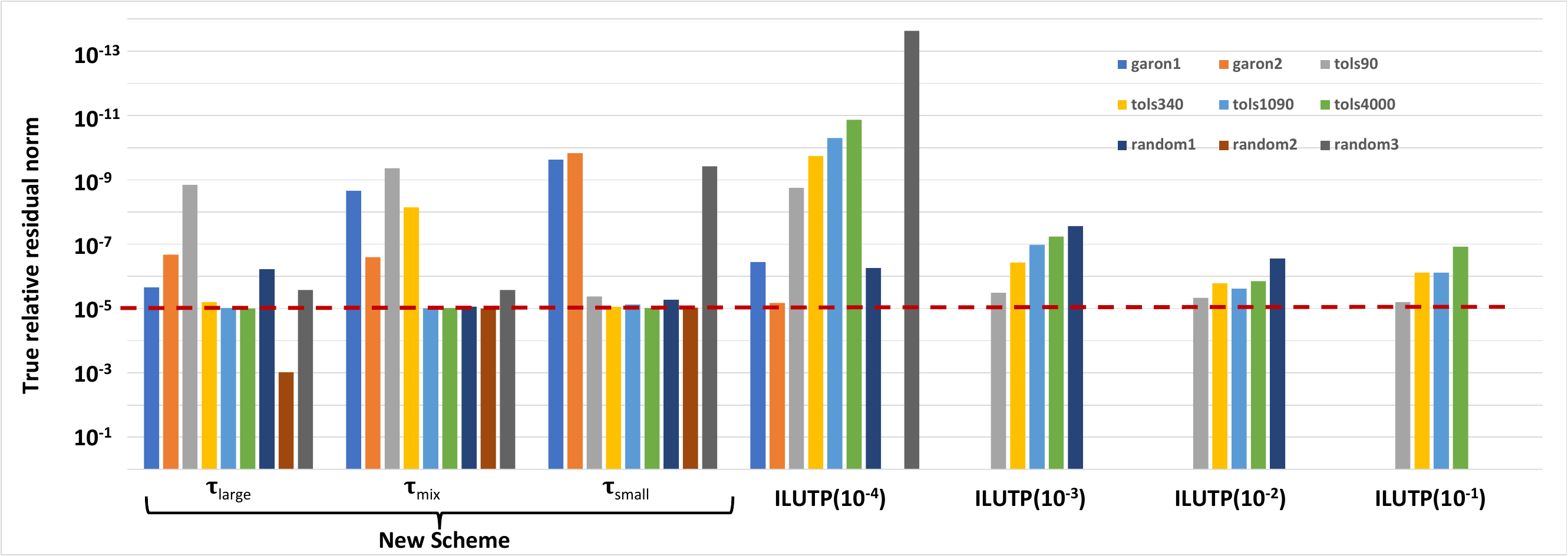}
  \caption{Final true relative residual norms of the new scheme and the baseline for the structurally symmetric case}
  \label{fig:truerelativeresidual:general}
\end{figure}

\begin{sidewaystable}[htbp]\centering
\begin{tabular}{r|lll|llll}
&\multicolumn{3}{c}{New Scheme} &\multicolumn{4}{c}{ILUTP($\tau$)} \\
\# & $\tau_{large}$ & $\tau_{mix}$ & $\tau_{small}$  & $10^{-4}$ & $10^{-3}$ &$10^{-2}$  & $10^{-1}$ \\ \hline
1 & 1,252,275 & 1,252,275 & 2,204,857  & {\bf 623,798} &$\dagger$ &$\dagger$ &$\dagger$ \\
2 & 18,037,895 & 18,037,895 & 31,745,441  & {\bf 5,345,750} &$\dagger$ &$\dagger$ &$\dagger$ \\
3 & {\bf 162} & {\bf 162} & {\bf 162}  & 1,797 & 1,642  & 1,229  & 640  \\
4 & {\bf 612} & {\bf 612} & {\bf 612} & 2,497 & 2,342  & 1,860 &  1,205 \\
5 &{\bf 1,962}  & {\bf 1,962} &  {\bf 1,962}  & 4,597 & 4,442 & 3,618 &2,882  \\
6 & {\bf 7,200}  & {\bf 7,200}  &  {\bf 7,200}  & 12,745 & 11,603  & 10,233 & 9,344 \\
7 & {\bf 1,885} & 1,688 & 1,875  & 27,955 & 27,527  & 25,366 &$\dagger$ \\
8 & $114\ddagger$  & {\bf 90} & 134  & $\dagger$ & $\dagger$  & $\dagger$ &$\dagger$ \\
9 & {\bf 320}  & {\bf 320} & {\bf 320} & 5,298 & $\dagger$  & $\dagger$ &$\dagger$ \\
\end{tabular}
\caption{Number of nonzeros of the preconditioners for the general case}
\label{tab:general:results:memory}
\end{sidewaystable}
In Table~\ref{tab:general:results:inner}, we present the average number of inner fGMRES(10), LSQR and MRS iterations for the new scheme. For all problems (except garon\_1 and garon\_2) the average number of inner fGMRES(10) iterations is small, even if $\tau_{large}$ and $\tau_{mix}$ are used. As in the other problem classes, the average number of LSQR iterations can be quite large for the set of general problems as well. The average number of MRS iterations is reasonably small for all $\tau_{\{large, mix, small\}}$.
\begin{table}[htbp]\centering
\begin{tabular}{r|lll|lll|lll}
\# &\multicolumn{3}{c}{$\tau_{large}$} &\multicolumn{3}{c}{$\tau_{mix}$}  & \multicolumn{3}{c}{$\tau_{small}$}  \\
 & fgmres & lsqr & mrs & fgmres  & lsqr & mrs & fgmres  & lsqr & mrs \\ \hline
1 & 671.0 & 67.3  & 6.3 & 753.3 & 86.8  & 7.0 & 6.0 &107.3  & 7.0  \\
2 & 987.0 & 113.1  & 1.3 & $\ddagger$& 152.1  & 1.9 & 4.0 & 197.0 & 1.9  \\
3 & 1.0 & 11.3  & 1.0   & 1.0 & 11.8  & 1.0 & 1.0 & 12.0 & 1.0  \\
4 & 1.0  & 45.4  & 1.0  & 1.0 &  92.3 & 1.0 & 1.0 & 94.5 & 1.0  \\
5 & 1.0 & 130.1  &1.0  &1.0 & 392.0  & 1.0 & 1.0 & 500.5 & 1.0  \\
6 &1.0  & 216.5  & 1.0 &1.0 & 500.5  &1.0  & 1.0 & 500.5 & 1.0  \\
7 & 8.0 & 87.1  & 8.0  & 7.0 & 95.4  & 10.6 &  7.0&  88.0 & 9.7   \\
8 & 7.3 & 70.1  & 7.3  &7.0 & 155.7  & 7.5 & 6.2 & 373.9 & 8.3  \\
9 & 2.0 & 178.6  & 7.4 & 2.0& 527.0  &1.0  & 2.0 & 611.0 & 1.0
  \end{tabular}
  \caption{Average number of inner fGMRES,  LSQR and MRS iterations for the new scheme (general case)}
\label{tab:general:results:inner}
\end{table}
Figure~\ref{fig:truerelativeresidual:general} shows the final true relative residual norms of the new scheme and the baseline method for the general case. The  horizontal dashed line visualizes the stopping tolerance where the norm of the relative residual is equal to $10^{-5}$.

\section{Conclusions}
\label{sec:conclusions}

We have introduced a new class of multilayer iterative schemes for solving sparse
linear systems in saddle point structure. We have presented a  theoretical analysis for the  (approximate) nullspace method. We have demonstrated the effectiveness of the new schemes by solving linear systems that arise in a variety of applications and compare them to a classical preconditioned iterative scheme. The new schemes may benefit from further preconditioning in the  Krylov subspace methods for solving over/underdetermined linear least squares problems. This is currently under investigation. 
Furthermore, a parallel implementation of the new schemes for large-scale problems is also planned as  a future work.

\section*{Acknowledgements}
The first author was supported by the Alexander von Humboldt Foundation for a research stay at Technische Universität Berlin.

\end{document}